\documentclass[11pt,reqno]{amsart}
\usepackage{amscd,graphicx,psfrag}
\usepackage{amssymb}
\usepackage{color}

\renewcommand{\phi}{\varphi}

\renewcommand{\det}{\operatorname{det}}

\newtheorem*{Theorem}{Theorem}

\theoremstyle{definition}

\begin{document}

\title{On instanton homology of corks $W_n$}

\author{Eric Harper}
\address{Department of Mathematics, McMaster University,
Hamilton, Ontario}
\email{eharper@math.mcmaster.ca}

\thanks{The author was partially supported by a CIRGET postdoctoral fellowship and a postdoctoral fellowship from McMaster University.}


\begin{abstract}We consider a family of corks, denoted $W_n$, constructed by Akbulut and Yasui.  Each cork gives rise to an exotic structure on a smooth $4$-manifold via a twist $\tau$ on its boundary $\Sigma_n = \partial W_n$.  
We compute the instanton Floer homology of $\Sigma_n$ and show that the map induced on the instanton Floer homology by $\tau: \Sigma_n \rightarrow \Sigma_n$ is non-trivial.
\end{abstract}

\maketitle
 

\section{Introduction}

In \cite{AY}, Akbulut and Yasui defined
a cork $C$ as a compact Stein $4$-manifold with boundary together with an involution $\tau: \partial C \rightarrow \partial C$ which
extends as a self-homeomorphism of $C$ but not as a self-diffeomorphism.  In addition, $C \subset X$ is a cork of a smooth $4$-manifold $X$ if cutting $C$ out and regluing it via $\tau$ changes the diffeomorphism type of $X$.

\begin{figure} [h]
\hfill
\begin{minipage}[t]{.45\textwidth}
\begin{center}
\def\svgwidth{95pt}
\input{wn-framed-link.pdf_tex}

\caption{$W_n$}
\label{WN}
\end{center}
\end{minipage}
\hfill
\begin{minipage}[t]{.45\textwidth}
\begin{center}
\def\svgwidth{100pt}
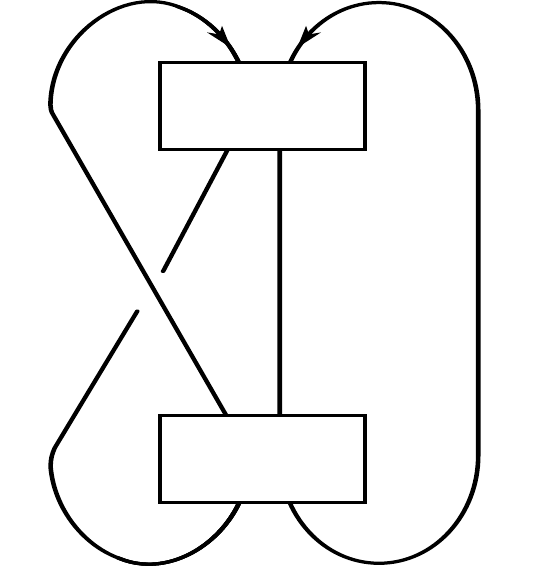
\caption{$\Sigma_n$}
\label{LN3}
\end{center}
\hfill
\end{minipage}
\hfill
\end{figure}

\begin{figure} [h]
\hfill
\begin{minipage}[t]{.45\textwidth}
\begin{center}
\def\svgwidth{115pt}
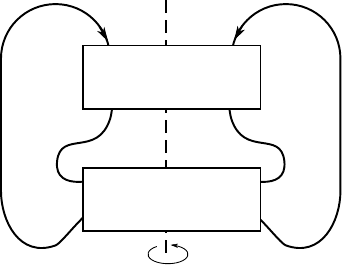
\caption{$L_n$}
\label{LN2}
\end{center}
\end{minipage}
\hfill
\begin{minipage}[t]{.45\textwidth}
\begin{center}
\def\svgwidth{125pt}
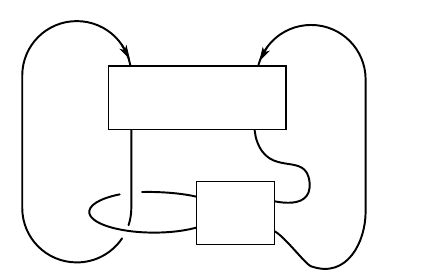
\caption{$S^3$}
\label{KN}
\end{center}
\hfill
\end{minipage}
\hfill
\end{figure}

We will consider the family of corks $W_n$, $n \geq 1$, obtained by surgery on the link in Figure \ref{WN} where a positive integer $m$ in a box indicates $m$ right-handed half-twists.  The boundary $\Sigma_n$ of $W_n$ is the integral homology $3$-sphere with surgery description as in Figure \ref{LN3}.  The involution $\tau : \Sigma_n \rightarrow \Sigma_n$ interchanges the two components of the link in Figure \ref{LN3}.  It is best seen when the underlying link $L_n$ is drawn symmetrically, as in Figure \ref{LN2}.  
Note that the quotient manifold $\Sigma'_n=\Sigma_n / \tau$ is homeomorphic to $S^3$ so $\Sigma_n$ can be viewed as a double branched cover of $S^3$ with branch set $k_n$ as shown in Figure \ref{KN}.

The goal of this paper is to study the instanton Floer homology $I_*(\Sigma_n)$ and the map $\tau_* : I_*(\Sigma_n) \rightarrow I_*(\Sigma_n)$ induced on it by $\tau$.  

%
%

\begin{Theorem}
$(1)$ For every integer $n \geq 1$, the instanton Floer homology group $I_j(\Sigma_n)$, $j=0,\ldots,7$, is trivial if $j$ is even, and is a free
abelian group of rank $n(n+1)(n+2)/6$ if $j$ is odd.

$(2)$
The homomorphism $\tau_* : I_*(\Sigma_n) \rightarrow I_*(\Sigma_n)$ is non-trivial for all $n\geq1$. 
\end{Theorem}

The first example of an involution acting non-trivially on the instanton Floer homology of an irreducible homology $3$-sphere was given in \cite{AK0} and \cite{Saveliev}; in fact, that example was exactly our $\tau: \Sigma_1 \rightarrow \Sigma_1$.
The technique we use to show non-triviality of $\tau_*$ is the same as the technique that was used in \cite{RS2} to reprove the result of \cite{Saveliev}: compare the Lefschetz number of $\tau_*: I_*(\Sigma_n) \rightarrow I_*(\Sigma_n)$ with the Lefschetz number of the identity map.  If the two are different, then the involution must be non-trivial.   
For any integral homology $3$-sphere $\Sigma$, the Lefschetz number of the identity equals the Euler characteristic of $I_*(\Sigma)$, which by Taubes \cite{Taubes} is twice the Casson invariant $\lambda(\Sigma)$.  Ruberman and Saveliev \cite{RSpro} showed that the Lefschetz number of $\tau_*$ equals twice the equivariant Casson invariant $\lambda^ \tau(\Sigma)$, defined in \cite{CS}.  Therefore the non-triviality of $\tau_*$ will follow as soon as we show that $\lambda(\Sigma_n) \neq \lambda^ \tau (\Sigma_n)$.  The calculation of $I_*(\Sigma_n) \rightarrow I_*(\Sigma_n)$ is done using surgery techniques.  

It should be noted that Akbulut and Karakurt proved a Heegaard Floer analogue of this result.  In \cite{AK} they showed that the involution $\tau: \Sigma_n \rightarrow \Sigma_n$ acts non-trivially on the Heegaard Floer homology group $HF^{+}(\Sigma_n)$.  

\subsection*{Acknowledgements}
I would like to thank Nikolai Saveliev for many helpful conversations and useful comments on this work, and \c{C}a\u{g}ri Karakurt for bringing
\cite{Maruyama} to our attention.

\section{The instanton Floer homology of $\Sigma_n$}\label{Splice}

Let $\Sigma$ be an oriented integral homology $3$-sphere. The instanton homology groups $I_*(\Sigma)$ are eight abelian groups arising as the Floer homology  of the Chern-Simons functional on the space of irreducible $SU(2)$ connections on $\Sigma$ modulo gauge equivalence.  The Euler characteristic of $I_*(\Sigma) $ is twice the Casson invariant $\lambda(\Sigma)$, see Taubes \cite{Taubes}.

Let $\Sigma_n(p)$ be the integral homology sphere with surgery description obtained by replacing the $0$-framing of the unknot on the left-hand side of Figure \ref{LN3} by a $p$-framing.  In \cite{Saveliev}, Saveliev used the Floer exact triangle to show  that the instanton homology groups of $I_*(\Sigma_1)(p)$ are independent of $p$.  The same argument holds for $I_*(\Sigma_n)(p)$.  

The homology $3$-sphere $\Sigma_n(2n+2)$ was shown by Maruyama \cite{Maruyama} to be homeomorphic to the Brieskorn homology sphere $\Sigma(2n+1, 2n+2, 2n+3)$.  In Theorem $10$ of \cite{Saveliev2}, Saveliev proved that the Floer homology group 
$I_j(\Sigma(2n+1, 2n+2, 2n+3))$ is trivial when $j$ is even, and is isomorphic to a free abelian group of rank $n(n+1)(n+2)/6$ when $j$ is odd.  This completes the calculation of $I_*(\Sigma_n)$.  In particular, the Casson invariant of $\Sigma_n$ is given by 
$$
\lambda(\Sigma_n) = -\frac{1}{3}n(n+1)(n+2). 
$$


\section{The Casson invariant}

Although the instanton Floer homology groups of $\Sigma_n$ are known, and therefore so is its Casson invariant, we can compute
$\lambda(\Sigma_n)$  using topological methods.

Let $L$ be a framed $2$-component link in $S^3$, and assume that the $3$-manifold $\Sigma$ resulting from surgery on $L$ is a homology $3$-sphere.  Boyer and Lines \cite{BL} compute $\lambda(\Sigma)$ as a sum of derivatives of the multivariable Alexander polynomial $\Delta_L$ of the underlying oriented link $L$ (and those of its sublinks), and Dedekind sums that only depend on the framings of the link components.  

In our case, both components of $L_n$ are framed by zero, and the Boyer-Lines \cite{BL} formula  for $\lambda(\Sigma_n)$ is simply 
\begin{equation}\label{Casson}
\lambda(\Sigma_n) =  -\frac{1}{\det(B)}\frac{\partial^2\Delta_{L_n}}{\partial x \partial y} (1,1), 
\end{equation}
where $B$ is the framing matrix for $\Sigma_n$.  Thus to compute $\lambda(\Sigma_n)$ we will need only to compute $\Delta_{L_n}$.

\subsection{The Alexander polynomial and Conway potential function}

Rather than computing $\Delta_L$ directly we will consider the related Conway potential function $\nabla_L$.
Given an oriented link $L$ in $S^3$, normalize $\Delta_L(x,y)$ using the Conway potential function $\nabla_L (x,y)$ of Hartley \cite{Hartley} by requiring that
\begin{equation}\label{DeltaNabla}
\Delta_L(x^2, y^2) = \nabla_L(x,y).
\end{equation}
Note that (\ref{DeltaNabla}) implies that
\begin{equation}\label{Casson2}
\frac{\partial^2\Delta_L}{\partial x \partial y} (1,1) = 
\frac{1}{4}\frac{\partial^2\nabla_L}{\partial x \partial y} (1,1),
\end{equation}  
hence we only need to compute the Conway potential function of $L$ and its partial derivatives.

\begin{figure}[h]
\centering
\def\svgwidth{150pt}
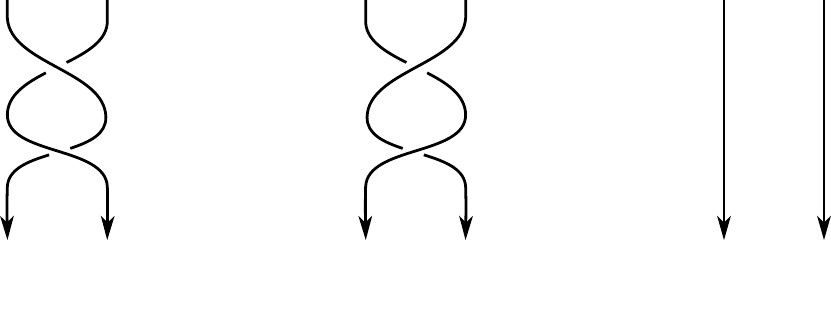
\caption{}
\label{RR}
\end{figure}

The Conway potential function $\nabla_L$ enjoys the replacement relation
\begin{equation}\label{RR1}
\nabla_{\ell} + \nabla_{r} = (xy + x^{-1} y^{-1})\nabla_{s},
\end{equation}
where $\nabla_{\ell}$, $\nabla_{r}$, and $\nabla_{s}$ are Conway functions of links that differ only in a neighborhood of a single crossing as shown in Figure \ref{RR}.  Note that the arcs may belong to the same component of $L$ or to different components.  If the three links differ as in Figure \ref{RR} but with one of the arcs oppositely oriented, then we have the relation
\begin{equation}\label{RR2}
\nabla_{\ell} + \nabla_{r} = (xy^{-1} + x^{-1} y)\nabla_{s}.
\end{equation}
Also, we will note that the Conway potential function vanishes for split links and is equal to $1$ for the right-handed Hopf link.

In order to calculate the Conway potential function of $L_n$ we will use the replacement relations (\ref{RR1}) and (\ref{RR2}) to produce 
a linear recurrence which we will then solve.

\subsection{The recurrence relation}

\begin{figure}[h]
\centering
\def\svgwidth{100pt}
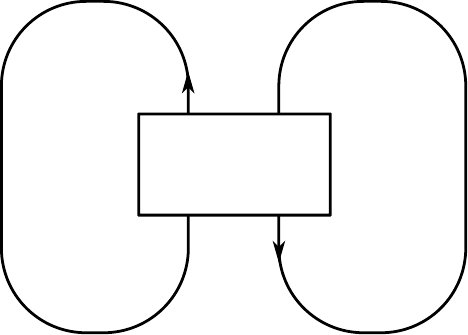
\caption{$H_n$}
\label{HN}
\end{figure}

Let $f_n=\nabla_{L_n}$ and $g_n = \nabla_{H_n}$,
where $H_n$ is the link in Figure \ref{HN}.  A straightforward calculation shows that 
\[
g_n=\left(\frac{1}{r_2-r_1}\right)r_1^n+\left(\frac{1}{r_1-r_2}\right)r_2^n
\]
with $v=xy^{-1} + x^{-1}y$, $r_1=(v+\sqrt{v^2-4})/2$, and $r_2=(v-\sqrt{v^2-4})/2$.
Using Hartley's replacement relations (\ref{RR1}), we change crossings of $L_n$ two at a time,
until we have undone the upper tangle in Figure \ref{LN3}.  We obtain the recurrence  
$f_{n+2} = -f_n+ uf_{n+1}$ with initial conditions
$f_0=-g_{n+1} + ug_n$ and 
$f_1=-ug_{n+1}+(u^2-1)g_n$,
where $u = xy + x^{-1} y^{-1}$.

Solving the recurrence relation, we obtain a formula for $f_n$ in terms of $x,y,$ and $n$,

\[
f_n=\left(f_0+\frac{f_0s_1-f_1}{s_2-s_1}\right)s_1^n+\left(\frac{f_1-f_0s_1}{s_2-s_1}\right)s_2^n
\]
\\
\noindent
where $s_1=(u+\sqrt{u^2-4})/2$ and $s_2=(u-\sqrt{u^2-4})/2$.  
We then find an explicit formula for $\lambda(\Sigma_n)$ by having \emph{Maple} \cite{Maple} differentiate $f_n$ twice and putting together (\ref{Casson}) and (\ref{Casson2}).  The answer is
\begin{equation}\label{CassonSigmaN}
\lambda(\Sigma_n) = -\frac{1}{3}n(n+1)(n+2).
\end{equation}

\section{The Equivariant Casson invariant}

Let $\tau : \Sigma \rightarrow \Sigma$ be an orientation preserving involution on a homology $3$-sphere, and suppose that the fixed point set of $\tau$ is non-empty.  Then the quotient manifold $\Sigma' = \Sigma/\tau$ is a homology $3$-sphere, and the projection map $\Sigma \to \Sigma'$ is a double branched cover with branch set a knot $k\subset \Sigma'$.  In \cite{CS}, Collin and Saveliev computed the equivariant Casson invariant $\lambda^ \tau (\Sigma)$ in terms of $\lambda(\Sigma')$ and the knot signature $\sigma(k)$.
When $\Sigma' = S^3$, we have simply 
\[
\lambda^ \tau(\Sigma) = \frac{1}{8} \sigma(k).
\]

Since we know that $\lambda(\Sigma_n)$ is decreasing as $n \rightarrow \infty$, see (\ref{CassonSigmaN}), our strategy for showing that $\lambda^{\tau}(\Sigma_n) \neq \lambda(\Sigma_n)$ will be to show 
that $\sigma(k_n)/8$ is bounded from below by a function strictly greater than $\lambda(\Sigma_n)$, where $k_n$
is the knot shown in Figure \ref{KN}.

\subsection{Bounding knot signatures}

\begin{figure}[h]
\centering
\def\svgwidth{90pt}
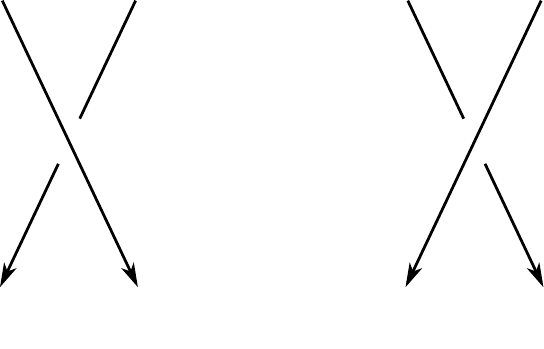
\caption{}
\label{Giller}
\end{figure}

The knot signature of a knot $k \subset S^3$ may be bounded from below using the formula
\begin{equation}\label{SigBound}
\sigma(k_r) \leq \sigma(k_{\ell}) \leq \sigma(k_r) + 2, 
\end{equation}
see Conway \cite{Conway} or Giller \cite{Giller}.  Here $k_r$ and $k_{\ell}$ are knots that only differ in a neighborhood of a crossing as shown in Figure \ref{Giller}. Note that our sign convention is opposite of Giller's.  
By (\ref{SigBound}), the signature of a knot $k$ is bounded from below by negative twice the number of right-handed crossings that must be changed in order to undo the knot $k$.  Note that we may need to change some left handed crossings while undoing $k$ but this will not contribute to our estimate.

\begin{figure}[h]
\centering
\def\svgwidth{225pt}
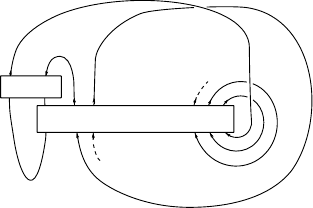
\caption{$k_n$}
\label{BD}
\end{figure}

We will now apply this observation to the knot $k_n$ shown in Figure \ref{KN}.  In order to see $k_n$ more clearly, we will isotope the link diagram in Figure \ref{KN} so that the $1$-framed curve is interchanged with the branch set $k_n$ and then blow down the $1$-framed curve.  The blow down has the effect of a full left-handed twist on the $2n+4$ strands passing through the $1$-framed curve, see Figure \ref{BD}.   

Note that when $n$ is odd, Figure \ref{BD} is a diagram of $k_n$ where no right-handed crossings must be changed in order to undo the knot, hence $\sigma(k_n)$ must be non-negative.  When $n$ is even, the left-most strand of the strands being twisted in Figure \ref{BD} is oriented oppositely of the other strands being twisted.  In this case there are $2n+3$ right-handed crossings that must be changed.  We change them to arrive at a knot with only left-handed crossings that further need to be changed.  Thus the signature of $k_n$ is bounded from below by $-4n-6$, and we have shown that
$$
\lambda(\Sigma_n)=-\frac{n(n+1)(n+2)}{3} < -\frac{4n+6}{8} \leq \lambda^{\tau}(\Sigma_n).
$$ 

Lastly, we remark that for small $n$ the equivariant Casson invariant can be computed explicitly using the following technique.  If we undo the upper left tangle in Figure \ref{BD} by changing crossings, then we will arrive at a knot that is isotopic to a torus knot.  If $n$ is odd, then the corresponding torus knot is the $T(2n+4, 2n+3)$ torus knot.  If $n$ is even, then the corresponding torus knot is $T(2n+2, 2n+1)$.   If $n \leq 4$, then we need to change at most $3$ crossings, and consequently the signature of $k_n$ differs from that of the corresponding torus knot by at most $\pm 6$.
  
For example, if $n = 1$ or $n=2$, then the corresponding torus knot is $T(6,5)$ and the signature $\sigma(T(6,5))$ is $16$.    Since the signature of $k_n$ is divisible by $8$, we have that $\sigma(k_1)=\sigma(k_2)=16$.  Similarly,  if $n=3$ or $n=4$, then the corresponding torus knot is 
$T(10,9)$ and 
$$
\sigma(k_3)=\sigma(k_4)=\sigma(T(10,9))=48.
$$

\end{document}